\documentclass[12pt]{article}
\usepackage{setspace}
\usepackage{textcomp}
\usepackage{amsmath, amsthm, amssymb}
\usepackage[left=1.5in,top=1in,right=1in,nohead]{geometry}
\usepackage{hyperref}
\usepackage{latexsym}
\usepackage{mathrsfs}
\usepackage{rawfonts}
\usepackage[dvips]{graphicx}
\usepackage[all]{xy}
\usepackage{stmaryrd}

\usepackage[sc,osf]{mathpazo}

\def\bct{\begin{center}}
\def\ect{\end{center}}
\def\beg{\begin}

\def\<{\langle}
\def\>{\rangle}

\def\mbb{\mathbb}
\def\mbbz{\mathbb Z}

\def\ni{\noindent}

\def\tn{\textnormal}

\newtheorem{thm}{Theorem}[section]
\newtheorem{lem}[thm]{Lemma}

\title{Complex and Symplectic Structures on Panelled Web 4-Manifolds} \author{H\"ulya Arg\"uz \and Mustafa Kalafat}

\begin{document}
\maketitle
\begin{abstract} We analyze the symplectic and complex structures on the panelled web 4-manifolds of \cite{handle}. In particular, we give infinite families of examples of almost complex but not symplectic and not complex 4-manifolds in the non-simply connected case.
  \end{abstract}

\section{Introduction}

Panelled web $4$-manifolds are introduced in \cite{handle} which have a natural locally conformally flat (LCF) structure with negative scalar curvature. The underlying smooth manifolds have the nice property that they can be used to produce self-dual metrics with small signatures. See \cite{ako} for further information. In this paper we analyze complex and symplectic structures on these type of manifolds. Summarizing the Theorems \ref{acx1}, \ref{acx2}, \ref{acx3} and \ref{acx4} we obtain the following result.

\vspace{.05in}

\noindent {\bf Theorem A.}
{\em The manifolds $M^1_g$, $M^2_{g,n}$, $M^3_{g,n}$ for $g,n>0$ and $M^4_n$ for $n>1$ are almost complex.}\footnote{Throughout the paper we have the assumption $g > 0$. Panelled 4-manifolds are built upon lower dimensional hyperbolic manifolds. These hyperbolic manifolds are produced by some special Kleinian groups which results in this condition. See 
\cite{maskitpwg} for more information.
} 

\vspace{.05in}

\ni Furthermore as an outcome of the search for a complex structure we get,

\vspace{.05in}

\noindent {\bf Theorem \ref{nocx}.}
{\em The manifolds $M^1_g$, $M^2_{g,n}$, $M^3_{g,n}$, $M^4_n$ do not admit any complex structure for all $g,n>0$.}

\vspace{.05in}

\ni Finally, the search for a symplectic structure yields the following.

\vspace{.05in}

\noindent {\bf Theorem \ref{nosymp}.}
{\em The manifolds $M^1_g$, $M^2_{g,n}$, $M^3_{g,n}$, $M^4_n$ do not admit any symplectic structure for all $g,n>0$.}

\vspace{.05in}

It follows from the above theorems that the non-simply connected four-manifolds we investigated all carry almost complex structures, but they are neither symplectic nor complex. 
The first kind of example of almost complex but noncomplex four-manifolds was constructed by A. Van de Ven. In his paper \cite{vandeven} he proves that for a four manifold $M$, if the Chern numbers $p=c_1^2[M]$ and $q=c_2[M]$ belong to a restricted domain then there does not exist a complex structure. Combining this with Milnor's result \cite{milnorcobordism}, namely if $p+q \equiv 0 \, (\tn{mod}\, 12)$ 
then an almost complex compact manifold can be constructed with these $p,q$; he constructs almost complex but noncomplex manifolds for special $p$ and $q$'s which are nonsimply connected and reducible, i.e. of the form of a connected sum.  Later on A. Howard \cite{alanhoward} gave examples using the same method of simply connected almost complex but noncomplex four manifolds. Shing-Tung Yau produces the first parallelizable manifolds of this kind in \cite{yau76} using Massey products, N. Brotherton provides more examples of parallelizable noncomplex almost complex manifolds in \cite{neilbrotherton} using the same techniques. These are basically constructed by taking the product of a 3-manifold with a circle. They all depend on the fact that a compact complex four manifold with even first Betti number is K\"ahler, and the real Massey products on a compact K\"ahler manifold vanish. So they construct compact manifolds with first even Betti numbers and nonvanishing Massey products, so that there is no way for them to be complex. Here we also provide examples of four-manifolds of this kind, namely almost complex, noncomplex and nonsymplectic manifolds. Ours have strictly negative (hence nontrivial) Euler characteristics, so that they are not parallelizable. In higher dimensions the existence of almost complex, noncomplex manifolds remains still as a conjecture. The reader may want to check \cite{fernandezgotaygray} for symplectic, noncomplex examples.

\vspace{.05in}

In Sections \S\ref{secacx} and \S\ref{seccx} we analyze the 
complex structures and in Section \S\ref{secsymp} we investigate symplectic structures on the manifolds.

\vspace{.05in}

{\bf Acknowledgements.} We would like to thank Anar Akhmedov for useful discussions. Also thanks to the referee for useful comments which improved the presentation. The figure is constructed by the IPE software of Otfried Cheong.

\section{Almost Complex Structures}\label{secacx}

The following result is a guideline for us to understand almost complex structures on 4-manifolds, see \cite{gs} for references.
\begin{lem}[\cite{hirzebruchhopf,wu}] \label{hhwu}
For a given 4-manifold $X$ with signature $\tau(X)$ and Euler characteristic $\chi(X)$, for any element $h \in H_2(X;\mbbz)$ satisfying the equation $h^2=3\tau(X)+2\chi(X)$ and the congruence $h \equiv w_2 \, (\tn{mod} \, 2)$ there is an almost complex structure $J$ on $TX$ with $h=c_1(X;J)$.
\end{lem}

\vspace{.05in}

Also we will be using the following fact, see \cite{gs} p.187 for details.

\beg{lem}\label{spinnability} Let $X$ be an oriented $4$-manifold (not necessarily compact).  If $H_1(X;\mbbz)$ has no $2$-torsion then $X$ is spin iff $Q_X$ is an even form. \end{lem}

\begin{figure}
\bct
\includegraphics[width=\textwidth]{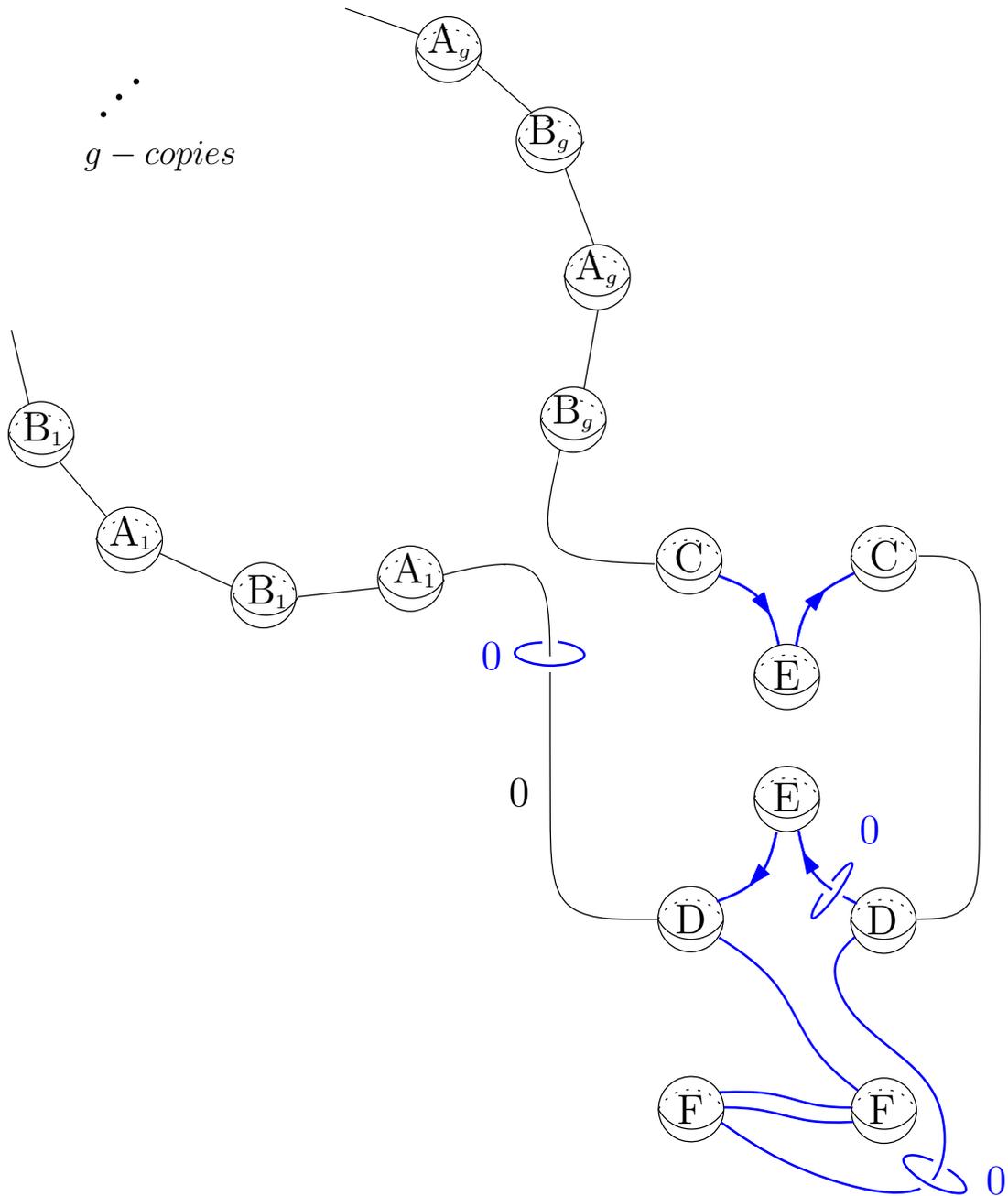}\\ \ect
  \caption{{\em {\small The Panelled Web 4-Manifolds $M^1_g$.} } } \label{sequenceofmetrics1}
\end{figure}
\vspace{.05in}

The first set of manifolds we want to work with are denoted by $M^1_g$, 
and can be seen through the Figure \ref{sequenceofmetrics1}. This is an instance of a {\em handlebody diagram.}
One can interpret the diagram as follows. Start with a blank page which stands for a copy of $\mbb R^3$, a part of $S^3$ the three dimensional sphere. The three-sphere bounds a $4$ dimensional ball from one side. 
As you can notice, the spheres in the picture appears in pairs. Through each pair we attach a $4$-dimensional $1$-handle, i.e. a copy of $D^1\times D^3$, where $D^k$ is a $k$-dimensional (closed) ball. After this process, we attach $4$-dimensional $2$-handles, i.e. $D^2\times D^2$ through the curves in the picture, each of which is a copy of a circle. Notice that a boundary piece of the $2$-handle is $S^1\times D^2$. After these two types of attachments, we need to attach the $3$ and $4$ handles, but these are uniquely attached so we do not mention them. In \cite{handle} the reader can find out the motivation to work with these manifolds. There, LCF metrics on these manifolds are considered. Here, we are only interested in the underlying smooth manifolds. Alternatively, for construction one begins with $\Sigma_{g,2}$, the twice punctured genus g surface, then cross it with the interval $I=[0,1]$  and then glue the boundary cylinders with each other with a flip. Finally make a so-called complex twist.  
This way one obtains the panelled web 3-manifold, and then by crossing this with $S^1$ and identifying its boundary one obtains the panelled web 4-manifold. A curious reader should check the reference to understand the relationship between the two constructions. 
These manifolds have the following characteristics.
$$\pi_1(M^1_g) = \langle a_1,b_1,...,a_g,b_g,c,d,e~|~a_1^{-1}b_1^{-1}a_1b_1\cdots a_g^{-1}b_g^{-1}a_gb_g cd^{-1}=ede^{-1}c^{-1}=df^3=1\rangle,$$
$$\tau(M^1_g)=0,~~ Q_{M^1_g}=H,~~ \chi(M^1_g)=-4g,$$
$$H_1(M^1_g;\mbbz)=\mbbz^{2g+2},$$ 
$$H_2(M^1_g;\mbbz) = \mbbz \oplus \mbbz.$$
Here, the fundamental group is generated by $1$-handles, and lowercase letters are used to denote the coressponding generators. 
As an application of the Lemma \ref{spinnability} they are all spin since they have even intersection form and there is no 2-torsion in $H_1(M^1_g;\mbbz)$. We will see that the Lemma is also applicable for the rest of the manifolds we will consider, so that they will be spin as well.
Our first result is the following.

\beg{thm} 
 \label{acx1} The manifolds $M^1_g$ are almost complex for all $g>0$. \end{thm}

\beg{proof} 
We know $H_2(M^1_g;\mbbz) = \mbbz \oplus \mbbz$. Say $h_1$ and $h_2$ are generators of $H_2(M^1_g;\mbbz)$. Let $(a,b)=ah_1+bh_2$ be an element of $H_2(M^1_g;\mbbz)$. If $(a,b)$ satisfies the above hypothesis we must have  
$(a,b)\cup (a,b)=3\tau(M^1_g)+2\chi(M^1_g)=3(0)+2(-4g)=-8g.$ 
But since we know from the intersection form of $M^1_g$, $h_1^2=h_2^2 =0$ and $h_1 h_2 =1$ we get $(a,b) \cup (a,b) = a^2 h_1 \cup h_1+b^2 h_2 \cup h_2+2abh_1 \cup h_2 = 2ab$. So we must have $2ab=-8g$, i.e. $ab= -4g$. Then for any fixed $g$, let $a=-2g$, $b=2$ so that $h=-2gh_1+2h_2 \in H_2(M^1_g;\mbbz)$. This gives an almost complex structure on $TM^1_g$ for any $g$. Note that the second condition in Wu's theorem, i.e. $h \equiv w_2 \, (\tn{mod}\, 2)$ is automatically satisfied for this choice of $h=(-2g,2)$, since we know that $M^1_g$ is spin hence $w_2 = 0$ and we have $h=-2gh_1+2h_2=0 \, (\tn{mod}\, 2)$. \end{proof}

\vspace{.05in}

In fact, the number of almost complex structures (upto homotopy) on each manifold of this family and upcoming families are finite. This is the consequence of that the Lemma \ref{hhwu} is actually an if and only if statement, i.e. almost complex structures are in correspondence with (or parametrized by) the elements of $H_2(M;\mbbz)$ satisfying the two conditions. One way to see this is that an almost complex structure turns the tangent bundle into a complex $2$-bundle, and $U(2)$ bundles are determined by the first two Chern classes, the latter of which is equal to the Euler class that is fixed for the tangent bundle so only $c_1$ conditions are important. See pp. $29-31$ of \cite{gs} for further information. For our family e.g. if $g$ is a prime number, then 
$(\pm 2g,\mp 2)$ and $(\pm 2,\mp 2g)$ are the only possibilities for $h$. In general $h^2$ conditions, which is $ab=-4g$ for this family restricts the possibilities to the finite case and depends on the number of prime factors of $g$.

\vspace{.05in}

The second family of manifolds we want to work on are the manifolds $M^2_{g,n}$. See the Figure 27 in \cite{handle} to see the explicit handlebody diagram. As one compares to the first family, the number of $CDE$ components is increased and the complex twist handle (the $F$ component) in the handlebody diagram for $M^1_g$ is omitted. The fundamental group and the other topological invariants are computed as follows. 
$$\pi_1(M^2_{g,n}) = 
\langle a_1,b_1,...,a_g,b_g,c_1,d_1,e_1,...,c_n,d_n,e_n~| \hspace{8cm}$$
$$\hspace{4.2cm}
~a_1^{-1}b_1^{-1}a_1b_1\cdots a_g^{-1}b_g^{-1}a_gb_gc_1\cdots c_n d_n^{-1}\cdots d_1^{-1} = e_i d_i e_i^{-1}c_i^{-1}=1 
\rangle.$$
$$\tau(M^2_{g,n})=0,~~ Q_{M^2_{g,n}}=H,~~ \chi(M^2_{g,n})=4-4g-4n.$$
$$H_1(M^2_{g,n};\mbbz) = \mbbz^{2g+2n}.$$
$$H_2(M^2_{g,n};\mbbz) = \mbbz \oplus \mbbz.$$

\ni From this set of constructions we obtain the following.

\beg{thm}\label{acx2} The manifolds $M^2_{g,n}$ are almost complex for all $g,n>0$. \end{thm} 

\beg{proof}We know $H_2(M^2_{g,n};\mbbz) = \mbbz \oplus \mbbz$. Say $h_1$ and $h_2$ are generators of $H_2(M^2_{g,n};\mbbz)$. Let $(a,b)=ah_1+bh_2$ be an element of $H_2(M^2_{g,n};\mbbz)$. If $(a,b)$ satisfies the above hypothesis we must have  
$(a,b)\cup (a,b)=3\tau(M^2_{g,n})+2\chi(M^2_{g,n})=3(0)+2(4-4g-4n)=8-8g-8n.$ 
But since we know from the intersection form of the manifold, $h_1^2=h_2^2 =0$ and $h_1 h_2 =1$ we get $(a,b) \cup (a,b) = a^2h_1 \cup h_1+b^2h_2 \cup h_2+2abh_1 \cup h_2 = 2ab$. So we must have $2ab=8-8g-8n$, i.e. $ab= 4-4g-4n$. Then for any fixed $g$, let $a=2-2g-2n$, $b=2$ so that $h=(2-2g-2n)h_1+2h_2 \in H_2(M^2_{g,n};\mbbz)$. This gives an almost complex structure on $TM^2_{g,n}$ for any $g$. Note that the second condition in Wu's theorem, i.e. $h \equiv w_2 \, (\tn{mod}\, 2)$ is automatically satisfied for this choice of $h=(2-2g-2n,2)$, since we know that $M^2_{g,n}$ is spin hence $w_2 = 0$ and we have $h=(2-2g-2n)h_1+2h_2=0 \, (\tn{mod}\, 2)$.
\end{proof}

\vspace{.05in}

In the third sequence of manifolds $M^3_{g,n}$ in \cite{handle} seen in Figure 29, 
$CE$-components are taken away and some new building blocks, namely trivial I-bundles over punctured annuli $\Sigma_{0,3}$ are attached through the handles $D_i$. 
These manifolds have,
$$\hspace{-1.3cm}\pi_1(M^3_{g,n}) = \langle a_1,b_1,...,a_g,b_g,k_1,l_1,m_1...,k_n,l_n,m_n~|~a_1^{-1}b_1^{-1}a_1b_1\cdots a_g^{-1}b_g^{-1}a_gb_gc_1\cdots c_n d_n^{-1}\cdots d_1^{-1}=1 \rangle.$$
$$\tau(M^3_{g,n})=0,~~ Q_{M^3_{g,n}}=(n+1)H,~~ \chi(M^3_{g,n})=4-4g-4n.$$
$$H_1(M^3_{g,n};\mbbz) = \mbbz^{2g+3n}.$$
$$H_2(M^3_{g,n};\mbbz) = \mbbz^{2+2n}.$$

\vspace{.05in}

\noindent Our next result is the following.

\beg{thm} 
 \label{acx3} The manifolds $M^3_{g,n}$ are almost complex for all $g,n>0$. \end{thm}   
 
\beg{proof} This time we have $H_2(M^3_{g,n};\mbbz) = \mbbz^{2+2n}$,\, 
so that $H_2(M^3_{g,n})$ is generated by $h_1,...,h_{2+2n}$
with $h_i \cup h_i=0$ for any $i$, and $h_i \cup h_{i+1}=h_{i+1} \cup h_i=1$ for odd $i$'s. So for instance we may take $h=2h_1+2h_2+2h_3-2gh_4+2h_5-2nh_6$. With this choice of $h$, one can obtain $h \cup h=8-8g-8n=3\tau(M^3_{g,n})+2\chi(M^3_{g,n})$. Again since $h$ is arranged to be even its mod $2$ reduction is zero, which is the Stiefel-Whitney class of the spin manifold. \end{proof}
\vspace{.05in}

In the last set of constructions of Panelled Web manifolds $M^4_n$ in \cite{handle} seen in Figure 30, many copies of
the new building blocks are attached to each other as a chain. By this way it is found that,
$$\pi_1(M^4_n) = \langle \cdots g_i, h_i, j_i, k_i, l_i, \cdots ,m ~|~ 
k_ih_ik_i^{-1}g_i = l_i^{-1} j_i l_i h_i = g_ih_ij_i=1 ~ \tn{for all} ~ i \leq n \rangle.$$
$$\tau(M^4_{n})=0,~~ Q_{M^4_{n}}=nH,~~ \chi(M^4_{n})=-2n.$$
$$H_1(M^4_n;\mbbz) = \mbbz^{2n+1}.$$
$$H_2(M^4_n;\mbbz) = \mbbz^{2n}.$$
\vspace{.05in}

\noindent Our next result is the following.

\beg{thm}\label{acx4} The manifolds $M^4_n$ are almost complex for all $n>1$ and do not admit any almost complex structure for $n=1$. \end{thm}
 
\beg{proof} We split the proof into the following cases.

\beg{itemize}

\item[Case 1:] Assume $n$ is even;
We have $H_2(M^4_{n};\mbbz) = \mbbz^{2n}$,\, 
so that $H_2(M^4_{n})$ is generated by $h_1,...,h_{2n}$
with $h_i \cup h_j=0$ for $j \neq i+1$ and $h_i \cup h_{i+1}=h_{i+1} \cup h_i=1$ for odd $i$'s. So for $n$ even we may take $h=2h_1+2h_2+2h_3-4h_4+2h_5+2h_6+2h_7-4h_8+...+2h_{2n-3}+2h_{2n-2}+2h_{2n-1}-4h_{2n}$ so that the coefficients take the values $(2,2,2,-4)$ repeatingly. With this choice of $h$, one can obtain $h \cup h=-4n=3\tau(M^4_{n})+2\chi(M^4_{n})$ and $h~(\tn{mod}\, 2)= w_2=0$. So, the hypothesis of Wu's Theorem is satisfied.  

\item[Case 2:] Assume $n$ is odd, $n \neq1$ and $n=2k+1$ with $k$ even (say $k=2t$, $t\in \mathbb{Z}^+$);
Again say $H_2(M^4_{n})$ is generated by $h_1,...,h_{2n}$ observe that $2n\geq10$ since with this choice of $n$, $n\geq5$. This time let $h=h_1+...+h_{2n-6}
-2h_{2n-5}+2th_{2n-4}-2h_{2n-3}+2th_{2n-2}-2h_{2n-1}+2th_{2n}$ so that the first $(2n-6)$ terms have coefficient $1$, and the coefficients for the last $6$ terms appear as $(-2,2t,-2,2t,-2,2t)$. Then $h \cup h=1+1+....+1-4t-4t-4t-4t-4t-4t=2n-6-24t=2n-6-3(8t)=2n-6-3(4k)=2n-6-3(2n-2)=-4n
=3\tau(M^4_{n})+2\chi(M^4_{n})$. Also, we know for all $i$, $h_i \neq 0=w_2~(\tn{mod}\, 2)$ since $h_i \cup h_{i+1}=1~(\tn{mod}\, 2)$ for odd $i$, but summing up $(2n-6)$ of them and adding the last terms with coefficients a multiple of $2$, guarantees the second condition in Wu's theorem. 

\item[Case 3:] Assume $n$ is odd, $n \neq1$ and $n=2k+1$ with $k$ odd (say $k=2t+1$, $t\in \mathbb{Z}$);
Let $h=h_1+...+h_{2n-2}-2h_{2n-1}+(6t+4)h_{2n}$ so that the first $(2n-2)$ terms have coefficient 1, and the coefficients for the last 2 terms are $(-2,6t+4)$. Then $h \cup h=1+1+....+1-12t-8-12t-8=2n-2-24t-16=2n-2-12(2t+1)-4=2n-2-12k-4=2n-2-6(2k+1)+2=2n-2-6n+2 =-4n=3\tau(M^4_{n})+2\chi(M^4_{n})$. Again we get $h=0~(\tn{mod}\, 2)$. 
\end{itemize}

This completes the proof for all $n$ such that $n\neq1$. In fact for $n=1$, $M^4_{n}$ are not almost complex since there can not exist any generator of $H_2(M^4_{n})$ satisfying the conditions of the Wu's Theorem in that case. \end{proof}


\section{Symplectic Structures}\label{secsymp}

Next we will investigate the symplectic structures on these manifolds. We follow the nice survey \cite{li}. There is a classification scheme for symplectic 4-manifolds similar to that of the Kodaira classification for complex surfaces. In the symplectic case, first of all one needs to find a substitute for the canonical bundle of the complex surface. This is easily generalized. Since associated to a symplectic manifold $(M,\omega)$, there is the contractible space of compatible almost complex structures. Using any compatible almost complex structure $J$ we define the symplectic Chern classes $c_k(M,\omega):=c_k(M,J)$ and hence the {\em symplectic
canonical class}  $-c_1(M,\omega)$ denoted sometimes by $-c_1(\omega)$ or $K_\omega$ which lies in $H^2(M;\mbbz)$. Next 
for a minimal symplectic 4-manifold $(M,\omega)$ 
we define the {\em Kodaira  dimension} in the following way. See \cite{li0,mss}.
$$\kappa(M,\omega) :=  \left\{ 
\begin{array}{cc}
-\infty& ~\tn{if}~\, c_1(\omega)\cdot[\omega]>0 ~~~\tn{or}~~~  c_1^2(\omega)<0\\
0      & ~\tn{if}~\, c_1(\omega)\cdot[\omega]=0  ~~\tn{and}~~  c_1^2(\omega)=0\\
1      & ~\tn{if}~\, c_1(\omega)\cdot[\omega]<0  ~~\tn{and}~~  c_1^2(\omega)^2=0\\
2      & ~\tn{if}~\, c_1(\omega)\cdot[\omega]<0  ~~\tn{and}~~  c_1^2(\omega)^2>0
\end{array} \right.$$

We will also make use of the following theorem. 
\beg{thm}[\cite{li}]Let M be  a  closed  oriented  smooth 4-manifold  and $\omega$ a  symplectic form on M  compatible with the orientation of M .  If $(M,\omega)$ is symplectically 
minimal, then 
\beg{enumerate}
\item The Kodaira dimension of $(M,\omega)$ is well-defined.
\item $(M,\omega)$  has Kodaira  dimension  $-\infty$ if and only if it is rational or ruled.
\item $(M,\omega)$   has Kodaira  dimension  0 if and only if $c_1(\omega)$  is a torsion class.  \end{enumerate}  \end{thm}

Notice that since all the manifolds $M^1_g, M^2_{g,n}, M^3_{g,n}, M^4_n$ have even intersection forms they can not contain a sphere of self-intersection $-1$, hence they are all minimal. Now we are ready to prove the following.

\beg{thm}\label{nosymp} None of the manifolds investigated above have  symplectic structure for $g,n>0$. \end{thm}

\beg{proof} Assume that the manifolds $M^1_g, M^2_{g,n}, M_{g,n}^3, M_{n}^4$ are symplectic, so that $\omega$ is a symplectic form on them. We have 
$$c_1^2(\omega)[M^1_g]     = (2\chi + 3\tau)(M^1_g)=-8g<0,$$ 
$$c_1^2(\omega)[M^2_{g,n}] = 8-8g-8n<0,$$ 
$$c_1^2(\omega)[M_{g,n}^3] = 8-8g-8n<0 ~~\tn{for}~~ g,n \geq 1 $$ and 
$$c_1^2(\omega)[M_{n}^4]   = 3\tau(M_{n}^4)+2\chi(M_{n}^4)=-2n<0 ~~\tn{for}~~ n \geq 1.$$ 
So all of the manifolds constructed must have Kodaira dimension $-\infty$, hence by the above theorem they must be rational or ruled if they are symplectic. If they were rational they should be $S^2\times S^2$ or $\mbb{CP}_2\sharp\,k\overline{\mbb{CP}}_2$ where k is any non-negative integer. This can not happen since the manifolds we consider  
are all non-simply connected. 
If they were ruled then they should be of the form $S^2 \times \Sigma_g \sharp\, k\overline{\mbb{CP}}_2$ where $\Sigma_g$ is the surface of genus $g$ and $k$ is any integer. This also fails to occur since the fundamental groups of $M^1_g, M^2_{g,n}, M_{g,n}^3, M_{n}^4$ which are stated in the previous section 
are not equal to the fundamental group of $S^2 \times \Sigma_g \sharp\, k\overline{\mbb{CP}}_2$. Therefore, they are neither rational nor ruled. So, they can not be symplectic.      
\end{proof}

\section{Complex Structures}\label{seccx}
In this section we will study the complex structures on our manifolds. The answer comes in the negative. For this study we will use the Enriques-Kodaira classification of surfaces. This classification is certainly an accumulation of works of many people, 
see \cite{bpv} or \cite{cx} for references.

\begin{lem}[Enriques-Kodaira Classification]  
Every compact, complex and connected surface has a minimal model in exactly one of the following classes $(1)$ to $(10)$. 

\vspace{.05in}

\begin{tabular}{|l|c|c|c|c|} 
\hline Class of $X$ & \tn{Kod} 
& $b_1$ 
&  $c^2_1$ & $c_2$ \\
\hline (1)Minimal rational surfaces & $-\infty$ & $0$ & $8$ or $9$ & $4$ or $3$ \\
\hline (2)Minimal surfaces of Kodaira's class VII & $-\infty$ &1&$\leq 0$ & $\geq 0$\\
\hline (3)Ruled surfaces of genus \ $\geq 0$ & $-\infty$ & $2g$ & $8(1-g)$ & $4(1-g)$\\
\hline (4)Enriques Surfaces    & $0$ &$0$ & $0$& $12$ \\
\hline (5)bi-elliptic surfaces & $0$ &$2$ & $0$&$0$ \\
\hline (6)Kodaira surfaces     & $0$ & $3 ~\tn{or}~1$ &$0$& $0$\\
\hline (7)K3 surfaces          & $0$ &$0$ & $0$&$24$ \\
\hline (8)Tori                 & $0$ &$4$&$0$ &$0$ \\
\hline (9)Minimal properly elliptic surfaces & $1$ & & $0$& $\geq 0$\\
\hline (10)Minimal surfaces of general type  & $2$ & even &$>0$ &$>0$ \\ 
\hline
\end{tabular}
\end{lem}



\vspace{.1in}

\noindent Next we will prove the following.
\beg{thm}\label{nocx} 
The manifolds $M^1_g$, $M^2_{g,n}$, $M^3_{g,n}$, $M^4_{n}$ do not admit any complex structure for all $g,n>0$. \end{thm}

 

\beg{proof} For all the manifolds $M^1_g, M^2_{g,n}, M_{g,n}^3, M_{n}^4$ we have already cited their characteristics in section \S\ref{secacx}, 
moreover we have $c_1^2(M^1_g)=-8g$, $c_1^2(M^2_{g,n})=8-8g-8n$, $c_1^2(M_{g,n}^3)=8-8g-8n$, $c_1^2(M_{n}^4)=-2n$ which are all negative for $n,g \geq 1$. According to the Kodaira-Enriques classification this
  can happen only when they are either in the second category i.e. a minimal surface of Kodaira's class VII or in the third category i.e. a ruled 
  surface
  $S^2 \times \Sigma_{g'} \sharp\, k\overline{\mbb{CP}}_2$ of genus $g' \geq 0$. The first case can 
  not happen since minimal surfaces of Kodaira's class VII have $b_1=1$ which is not the 
  case because $b_1(M^1_g)=2g+2\neq 1, b_1(M^2_{g,n})=2g+2n\neq 1, b_1(M^3_{g,n})=2g+3n\neq 1, b_1(M^3_n)=2n+1\neq 1$ for $n,g \geq 1$ (for the Betti numbers stated we refer to \cite{handle}). 
Finally, the ruled surface case is already ruled out as in the proof of Theorem \ref{nosymp}. 
Hence, the manifolds  
   $M^1_g, M^2_{g,n}, M_{g,n}^3, M_n^4$ are not complex. \end{proof}




{\small \beg{flushleft} 
\textsc{FB Mathematik, Universit\"at Hamburg, Bundesstrasse 55, 20146, Germany}\\
\textit{E-mail address:} \texttt{\textbf{huelya.arguez@\,math.uni-hamburg.de}}

{\small 
\beg{flushleft} \textsc{Tuncel\' \i  ~\" Un\' \i vers\' ites\' i, Turkia}\\
\textit{E-mail address:}  \texttt{\textbf{kalafg@\,gmail.com}} \end{flushleft}
}

\end{flushleft} }


\end{document}